\documentclass[12pt,reqno]{amsart}
\usepackage{enumerate, latexsym, amsmath, amsfonts, amssymb, amsthm, arydshln, color, mathrsfs}
\def\pmod #1{\ ({\rm{mod}}\ #1)}
\allowdisplaybreaks[4]
\def\Z{\Bbb Z}
\def\N{\Bbb N}

\def\bg{\bigg}
\def\({\bg(}
\def\){\bg)}
\def\t{\text}

\def\ord{{\rm ord}}

\theoremstyle{plain}
\newtheorem{theorem}{Theorem}

\newtheorem{lemma}{Lemma}

\theoremstyle{definition}
\newtheorem*{acknowledgment}{Acknowledgments}
\theoremstyle{remark}
\newtheorem{remark}{Remark}

\newcommand{\RNum}[1]

 \vspace{4mm}

\begin{document}

\hbox{Preprint}

\title
[{Numbers represented by restricted sums of four squares}]
{Numbers represented by restricted sums of four squares}

\author
[G.-L. Zhou and Y.-F. She]
{Guang-Liang Zhou and Yue-Feng She*}

\address{(Guang-Liang Zhou) Department of Mathematics, Nanjing
University, Nanjing 210093, People's Republic of China}
\email{guangliangzhou@126.com}

\address{(Yue-Feng She) Department of Mathematics, Nanjing
University, Nanjing 210093, People's Republic of China}
\email{she.math@smail.nju.edu.cn}

\thanks{*Corresponding author}

\begin{abstract} In this paper, we prove some results of restricted sums of four squares using arithmetic of quaternions in the ring of Lipschitz integers. For example, we show that every nonnegative integer $n$ can be written as $x^{2}+y^{2}+z^{2}+t^{2}$, where $x,y,z,t$ are integers and $x+y+2z+2t$ is a square (or a cube), which confirms a conjecture of Z.-W. Sun.
\end{abstract}
\maketitle

\section{Introduction}
\setcounter{lemma}{0}
\setcounter{theorem}{0}
\setcounter{corollary}{0}
\setcounter{remark}{0}
\setcounter{equation}{0}
In $1770$ Lagrange established his well-known four-square theorem, which shows that every positive integer $n$ can be  written as the sum of four squares. S. Ramanujan \cite{R} listed $55$ positive definite integral diagonal quaternary quadratic forms and claimed they can represent all positive integers. A few years later, L. E. Dickson \cite{D1} confirmed Ramanujan's assertion for $54$ of them and showed the remaining one, $x^{2}+2y^{2}+5z^{2}+5w^{2}$, can represent all positive integers other than $15$.

In \cite{S2}, Z.-W. Sun proposed the following interesting 1-3-5 conjecture, which is an extension of Lagrange's four-square theorem:

\medskip
\textit{
Every positive integer can be represented as $x^{2}+y^{2}+z^{2}+t^{2}$ with $x,y,z,t \in \N =\{0,1,2,\ldots\}$ such that $x+3y+5z$ is a square.
}
\medskip

\noindent Let $\Z$ denote the ring of integers. Y.-C. Sun and Z.-W. Sun \cite{SS} proved that any $n\in \N$ can be written as $x^{2}+y^{2}+z^{2}+t^{2}$ with $x,y,5z,5t\in \Z$ and $x+3y+5z$ a square using Euler's four-square identity. Later, H.-L. Wu and Z.-W. Sun \cite{W} investigated the integer version of 1-3-5 conjecture. Finally, with the help of  the Lipschitz integers, the 1-3-5 conjecture was completely proved by A. Machiavelo and N. Tsopanidis \cite{MT} recently.

Z.-W. Sun considered the more general problem, too.
Let $a,b,c,d$ be nonnegative integers and $\mathcal{S}$ be a subset of $\Z$.
For any positive integer $m$, consider
the system
\begin{equation}\label{eq1}
\begin{cases}
m=x^{2}+y^{2}+z^{2}+t^{2},\\
ax+by+cz+dt \in \mathcal{S}.\\
\end{cases}
\end{equation}
Z.-W. Sun \cite{S2,S3} proved some results concerning the solvability of the system \eqref{eq1} for the various $(a,b,c,d)$ and $\mathcal{S}$. For example, he \cite[Theorem 1 (iv)]{S3} proved that every positive integer can be represented  as $x^{2}+y^{2}+z^{2}+t^{2}$ with $x,y,z,t \in \Z$ and $x+y+z+t\in\{2^k;k\geq 1\}$. Furthermore, Sun \cite{S2,S3} also gave many related conjectures on the system \eqref{eq1}. For example, Sun's 24-conjecture \cite[Conjecture 4.7 (i)]{S3} asserts that every positive integer can be represented as $x^{2}+y^{2}+z^{2}+t^{2}$ with $x,y,z,t\in\N$ such that both $x$ and $x+24y$ are square.

Inspired by the above work, in the present paper, we have the following results.
\begin{theorem}\label{Th1.1}
Let $\mathcal{S}$ denote the set of all nonnegative cubic numbers. Then, for any nonnegative integer $m$, the system $(\ref{eq1})$ has an integer solution  if $(a,b,c,d)$ is among the nine quadruples
\begin{align*}&(1,1,2,2),(1,2,2,2),(2,2,3,0),(1,3,3,0),\\
&(1,2,4,0),(1,1,2,4),(2,3,4,0),(1,1,2,5),(1,2,3,5).
\end{align*}
\end{theorem}

\begin{remark}
	Z.-W. Sun \cite[Conjecture 4.4(iv)]{S2} first conjectured that each $n\in\N$ can be written as $x^{2}+y^{2}+z^{2}+t^{2}$ with $x,y,z,t\in \Z$ and $x+y+2z+2t$ a nonnegative cube. In addition, Z.-W. Sun \cite[Conjecture 4.4(iii)]{S3} conjectured that when $\mathcal{S}=\{0\}\cup\{\pm8^{k}: k\in \N\}$, the system (\ref{eq1}) has a natural solution if $(a,b,c,d)$ is among $(1,3,-3,0)$ and $(2,1,-1,-1)$.
\end{remark}

\begin{theorem}\label{Th1.2}
Let $\mathcal{S}=\{2^{k}:k\in\N\}$. Then, for any nonnegative integer $m$, the system $(\ref{eq1})$ has an integer solution if $(a,b,c,d)$ is among the five quadruples
$$(1,3,3,0),(1,2,4,0),(1,1,2,4),(1,1,2,5),(1,2,3,5).$$
\end{theorem}

\begin{remark}
Z.-W. Sun \cite[Conjecture 4.4(i) and Conjecture 4.4(ii)]{S3} conjectured that when $\mathcal{S}=\{2^{k}: k\in\N\}$, the system (\ref{eq1}) has a natural solution if $(a,b,c,d)$ is among
$$(1,3,-3,0),(4,-2,-1,0),(1,4,-2,0).$$

\end{remark}

\begin{theorem}\label{Th1.3}
Let $\mathcal{S}$ denote the set of all the squares. Then, for any nonnegative integer $m$, the system $(\ref{eq1})$ has an integer solution if  $(a,b,c,d)$ is among the four quadruples
$$(1,1,2,2),(1,2,2,2),(2,2,3,0),(2,3,4,0).$$
\end{theorem}

\begin{remark}
Z.-W. Sun \cite[Conjecture 4.3(iv) and Conjecture 4.12(i)]{S2} conjectured that when $\mathcal{S}$ denotes all the squares, the system (\ref{eq1}) has a natural solution if $(a,b,c,d)$ is among $(2,3,-2,0)$ and $(1,2,-1,-2)$. Moreover, Z.-W. Sun \cite[Conjecture 4.5(i)]{S3} conjectured that when $\mathcal{S}=\{0\}\cup\{4^{k}:k\in\N\}$, the system (\ref{eq1}) has a natural solution if $(a,b,c,d)=(2,-2,1,-1)$. In addition, Z.-W. Sun \cite[Conjecture 4.5(ii)]{S3} conjectured that when $\mathcal{S}=\{\pm4^{k}:k\in\N\}$, the system (\ref{eq1}) has a natural solution if $(a,b,c,d)$ is among $(1,1,2,-2)$ and $(2,3,-4,0)$.
\end{remark}

\begin{theorem}\label{Th1.4}
Suppose $m\in\N$.
\begin{itemize}
	\item  If $m\geq 3.74\times10^{9}$ and  $m\not\equiv0\pmod{16}$, $m$ can be written as $x^{2}+y^{2}+z^{2}+t^{2}$ with $x,y,z,t\in\N$ and $x+2y+3z+5t$ a square.
	\item If $m\geq 7.67\times10^{9}$ and  $m\not\equiv0\pmod{16}$, $m$ can be written as $x^{2}+y^{2}+z^{2}+t^{2}$ with $x,y,z,t\in\N$ and $2x+3y+4z$ a square.
\end{itemize}

\end{theorem}

\begin{remark}
Z.-W. Sun \cite[Conjecture 4.4(iii)]{S2} conjectured that each $m\in\N$ not of the form $7\times4^{2k+1}$ for any $k\in\N$ can be written as $x^{2}+y^{2}+z^{2}+t^{2}$ with $x,y,z,t\in\N$ and $x+2y+3z+5t$ a square. In addition, Z.-W. Sun\cite[Conjecture 4.14(ii)]{S4} conjectured that each $m\in\N$ not of the form $3\times2^{4k+1}$ for any $k\in\N$ can be written as $x^{2}+y^{2}+z^{2}+t^{2}$ with $x,y,z,t\in\N$ and $2x+3y+4z$ a square. To confirm these conjectures, we only need to check those positive integers less than the corresponding boundary.
\end{remark}

\section{some preparations}
\setcounter{lemma}{0}
\setcounter{theorem}{0}
\setcounter{corollary}{0}
\setcounter{remark}{0}
\setcounter{equation}{0}

We represent the system of equations$$
\begin{cases}
m=x^{2}+y^{2}+z^{2}+t^{2},\\
n=ax+by+cz+dt,\\
\end{cases}
$$
as ($a$-$b$-$c$-$d$) with the zeros left out (for example, we will use ($1$-$3$) instead of ($1$-$3$-$0$-$0$)).

For a quaternion $\alpha=a_{1}+a_{2}i+a_{3}j+a_{4}k$ ($a_{1},a_{2},a_{3},a_{4}\in\Z,\ i^{2}=j^{2}=k^{2}=ijk=-1$), let $N(\alpha)=a_{1}^{2}+a_{2}^{2}+a_{3}^{2}+a_{4}^{2}$ be the norm of $\alpha$ and $Re(\alpha)=a_{1}$ be the real part of $\alpha$.

Below is a useful tool we will need in our proof of main theorems.
\begin{lemma}\label{lem2.1}$($\cite[Theorem 1]{MT}$)$
Let $m,n,l \in \N$ be such that $n^{2}\leq lm$ and $lm-n^{2}$ is not of the form $4^{r}(8s+7)$ for any $r,s\in \N$. Then, for some $a,b,c,d\in \N$ satisfying that $a^{2}+b^{2}+c^{2}+d^{2}=l$, the system $($$a$-$b$-$c$-$d$$)$ has an integer solution.
\end{lemma}

 It is easy to verify that there are only two ways to write $10, 13, 17, 19,$ $ 21,  22, 29, 31$ and $39$ as sum of four squares. In fact,
\begin{align*}
10=&1^{2}+3^{2}=1^{2}+1^{2}+2^{2}+2^{2},\\
13=&2^{2}+3^{2}=1^{2}+2^{2}+2^{2}+2^{2},\\
17=&1^{2}+4^{2}=2^{2}+2^{2}+3^{2},\\
19=&1^{2}+1^{2}+1^{2}+4^{2}=1^{2}+3^{2}+3^{2},\\
21=&2^{2}+2^{2}+2^{2}+3^{2}=1^{2}+2^{2}+4^{2},\\
22=&2^{2}+3^{2}+3^{2}=1^{2}+1^{2}+2^{2}+4^{2},\\
29=&2^{2}+5^{2}=2^{2}+3^{2}+4^{2},\\
31=&2^{2}+3^{2}+3^{2}+3^{2}=1^{2}+1^{2}+2^{2}+5^{2},\\
39=&1^{2}+1^{2}+1^{2}+6^{2}=1^{2}+2^{2}+3^{2}+5^{2}.
\end{align*}

In light of Lemma \ref{lem2.1}, we would like to show that, under certain conditions, the existence of integer solutions of one system can induce the existence of integer solutions of the other system.

\section{Transformation of Solutions}
\setcounter{lemma}{0}
\setcounter{theorem}{0}
\setcounter{corollary}{0}
\setcounter{remark}{0}
\setcounter{equation}{0}

{\it Case} 1. ($2$-$3$-$3$)$\Rightarrow$($1$-$1$-$2$-$4$)

Let $\beta_{22}=2+3i+3j$ and $\gamma_{22}=x_{22}-y_{22}i-z_{22}j-t_{22}k$, where $(x_{22},y_{22},z_{22},t_{22})$ is an integer solution of the system ($2$-$3$-$3$).

One can easily check that
$$\beta_{22}(1+2i)=(1-2j)(-2-i-j-4k).$$
By computation, we have
\begin{footnotesize}
\begin{align*}
(1+2i)^{-1}\gamma_{22}(1-2j)=&\frac{x_{22}-2y_{22}-2z_{22}-4t_{22}}{5}+\frac{-2x_{22}-y_{22}+4z_{22}-2t_{22}}{5}i\\
&+\frac{-2x_{22}+4y_{22}-z_{22}-2t_{22}}{5}j+\frac{4x_{22}+2y_{22}+2z_{22}-t_{22}}{5}k.
\end{align*}
\end{footnotesize}

Note that $$N((1+2i)^{-1}\gamma_{22}(1-2j))=N(\gamma_{22})=m,$$ and $$Re((1-2j)^{-1}\beta_{22}(1+2i)(1+2i)^{-1}\gamma_{22}(1-2j))=Re(\beta_{22}\gamma_{22})=n.$$ Thus in this way we can construct a solution of the system ($1$-$1$-$2$-$4$) from a solution of ($2$-$3$-$3$) if
 $$ x_{22}-2y_{22}-2z_{22}+t_{22}\equiv0\pmod{5}.$$

Similarly, we have
\begin{align*}
\beta_{22}(1-2i)&=(1+2k)(4+i+j-2k),\\
\beta_{22}(1+2j)&=(1-2i)(-2-i-j+4k),\\
\beta_{22}(1-2j)&=(1-2k)(4+i+j+2k),\\
\beta_{22}(1+2k)&=(1+2i)(4+i+j+2k),\\
\beta_{22}(1-2k)&=(1+2j)(4+i+j-2k),
\end{align*}
and
\begin{footnotesize}
\begin{align*}
(1-2i)^{-1}\gamma_{22}(1+2k)=&\frac{x_{22}+2y_{22}+4z_{22}+2t_{22}}{5}+\frac{2x_{22}-y_{22}-2z_{22}+4t_{22}}{5}i\\
&+\frac{-4x_{22}+2y_{22}-z_{22}+2t_{22}}{5}j+\frac{2x_{22}+4y_{22}-2z_{22}-t_{22}}{5}k,\\
(1+2j)^{-1}\gamma_{22}(1-2i)=&\frac{x_{22}-2y_{22}-2z_{22}+4t_{22}}{5}+\frac{-2x_{22}-y_{22}+4z_{22}+2t_{22}}{5}i\\
&+\frac{-2x_{22}+4y_{22}-z_{22}+2t_{22}}{5}j+\frac{-4x_{22}-2y_{22}-2z_{22}-t_{22}}{5}k,\\
(1-2j)^{-1}\gamma_{22}(1-2k)=&\frac{x_{22}+4y_{22}+2z_{22}-2t_{22}}{5}+\frac{-4x_{22}-y_{22}+2z_{22}-2t_{22}}{5}i\\
&+\frac{2x_{22}-2y_{22}-z_{22}-4t_{22}}{5}j+\frac{-2x_{22}+2y_{22}-4z_{22}-t_{22}}{5}k,\\
(1+2k)^{-1}\gamma_{22}(1+2i)=&\frac{x_{22}+2y_{22}+4z_{22}-2t_{22}}{5}+\frac{2x_{22}-y_{22}-2z_{22}-4t_{22}}{5}i\\
&+\frac{-4x_{22}+2y_{22}-z_{22}-2t_{22}}{5}j+\frac{-2x_{22}-4y_{22}+2z_{22}-t_{22}}{5}k,\\
(1-2k)^{-1}\gamma_{22}(1+2j)=&\frac{x_{22}+4y_{22}+2z_{22}+2t_{22}}{5}+\frac{-4x_{22}-y_{22}+2z_{22}+2t_{22}}{5}i\\
&+\frac{2x_{22}-2y_{22}-z_{22}+4t_{22}}{5}j+\frac{2x_{22}-2y_{22}+4z_{22}-t_{22}}{5}k.
\end{align*}
\end{footnotesize}
Now, we can naturally get the following lemma.
\begin{lemma}\label{lem3.1}
 Given an integer solution $(x_{22}, y_{22},z_{22},t_{22})$ for the system $($$2$-$3$-$3$$)$,  if either of the following conditions holds:
\begin{itemize}
        \item $x_{22}-2y_{22}-2z_{22}+t_{22}\equiv0\pmod{5}$,
        \item $x_{22}+2y_{22}-z_{22}+2t_{22}\equiv0\pmod{5}$,
        \item $x_{22}-2y_{22}-2z_{22}-t_{22}\equiv0\pmod{5}$,
        \item $x_{22}-y_{22}+2z_{22}-2t_{22}\equiv0\pmod{5}$,
        \item $x_{22}+2y_{22}-z_{22}-2t_{22}\equiv0\pmod{5}$,
        \item $x_{22}-y_{22}+2z_{22}+2t_{22}\equiv0\pmod{5}$,
\end{itemize}
the system $($$1$-$1$-$2$-$4$$)$ has an integer solution.
\end{lemma}
In view of Lemma \ref{lem2.1}, we can induce the following lemma.

\begin{lemma}\label{lem3.2}
Let $m,n \in \N$ be such that $22m-n^{2}$ is nonnegative and not of the form $4^{r}(8s+7)$ for any $r,s \in \N$. If $22m-n^{2}\equiv0,1\ or\ 4 \pmod{5}$, the system $($$1$-$1$-$2$-$4$$)$ has an integer solution.
\end{lemma}

\medskip
\noindent{\it Proof}. By Lemma \ref{lem2.1}, either the system ($1$-$1$-$2$-$4$) or the system ($2$-$3$-$3$) has an integer solution. If the system ($1$-$1$-$2$-$4$) has an integer solution,  we are done. If not, we have an integer solution for the system  ($2$-$3$-$3$). Using the notations above, we have
\begin{align*}
\gamma_{22}\beta_{22}=n+A_{22}i+B_{22}j+C_{22}k,
\end{align*}
where
$$
\begin{cases}
n=2x_{22}+3y_{22}+3z_{22},\\
A_{22}=3x_{22}-2y_{22}+3t_{22},\\
B_{22}=3x_{22}-2z_{22}-3t_{22},\\
C_{22}=-3y_{22}+3z_{22}-2t_{22}.\\
\end{cases}
$$
By solving the equations and taking the coefficient modulus $5$, we can obtain that

\begin{small}
\begin{equation*}
\left(
  \begin{array}{c}  \vspace{1ex}
    x_{22}\\  \vspace{1ex}
    y_{22}\\  \vspace{1ex}
    z_{22}\\  \vspace{1ex}
    t_{22}\\
  \end{array}
\right)
\equiv
\left(
  \begin{array}{cccc} \vspace{1ex}
    1 & -1 & -1 & 0\\  \vspace{1ex}
    -1 & -1 & 0 & 1\\  \vspace{1ex}
    -1 & 0 & -1 & -1\\  \vspace{1ex}
    0 & -1 & 1 & -1
\end{array}
\right)
  \left(
  \begin{array}{c}  \vspace{1ex}
   n\\  \vspace{1ex}
   A_{22}\\  \vspace{1ex}
   B_{22}\\  \vspace{1ex}
   C_{22}\\
  \end{array}
\right)
\pmod{5}.
\end{equation*}
\end{small}
Hence
\begin{small}
\begin{equation*}
\left(
  \begin{array}{cccc}
    1 & -2 & -2 & 1\\
    1 & 2 & -1 & 2\\
    1 & -2 & -2 & -1\\
    1 & -1 & 2 & -2\\
    1 & 2 & -1 & -2\\
    1 & -1 & 2 & 2\\
    \end{array}
\right)
  \left(
  \begin{array}{c}  \vspace{1ex}
   x_{22}\\  \vspace{1ex}
   y_{22}\\  \vspace{1ex}
   z_{22}\\  \vspace{1ex}
   t_{22}\\
  \end{array}
\right)\equiv
\left(
  \begin{array}{cccc}
    0 & 0 & 2 & -1\\
    0 & 0 & 2 & 1\\
    0 & 2 & 0 & 1\\
    0 & 2 & 0 & -1\\
    0 & -1 & -2 & 0\\
    0 & -2 & -1 & 0\\
    \end{array}
\right)
  \left(
  \begin{array}{c}  \vspace{1ex}
   n\\  \vspace{1ex}
   A_{22}\\  \vspace{1ex}
   B_{22}\\  \vspace{1ex}
   C_{22}\\
  \end{array}
\right)
\pmod{5},
\end{equation*}
\end{small}
which means that if either of the following conditions holds:
\begin{itemize}
  \item $A_{22}\equiv \pm 2B_{22} \pmod{5}$,
  \item $A_{22}\equiv \pm 2C_{22} \pmod{5}$,
  \item $B_{22}\equiv \pm 2C_{22} \pmod{5}$,
\end{itemize}
the system ($1$-$1$-$2$-$4$) has an integer solution.

On the other hand, $A_{22}^{2}+B_{22}^{2}+C_{22}^{2}=22m-n^{2} \equiv 0,1 \ \t{or} \ 4\pmod{5}$. Note that
$$
\begin{cases}
0\equiv0+0+0\equiv0+1+4\pmod{5},\\
1\equiv0+0+1\equiv1+1+4\pmod{5},\\
4\equiv0+0+4\equiv1+4+4\pmod{5}.\\
\end{cases}
$$
One can easily see that one of $A^{2}+B^{2}$, $A^{2}+C^{2}$ and $B^{2}+C^{2}$ must be divided by $5$, hence one of the conditions must hold. This completes the proof. \qed
\medskip

{\it Case} 2. ($1$-$3$)$\Rightarrow$($1$-$1$-$2$-$2$) and ($2$-$3$)$\Rightarrow$($1$-$2$-$2$-$2$)

Let $\beta_{10}=1+3i$, $\beta_{13}=2+3i$, $\gamma_{10}=x_{10}-y_{10}i-z_{10}j-t_{10}k$ and $\gamma_{13}=x_{13}-y_{13}i-z_{13}j-t_{13}k$, where $(x_{10},y_{10},z_{10},t_{10})$ is an integer solution of the system ($1$-$3$) and $(x_{13},y_{13},z_{13},t_{13})$ is an integer solution of the system ($2$-$3$).

One can easily check that
\begin{align*}
\beta_{10}(1+i+j)=&(1+i+j)(1+i+2j+2k),\\
\beta_{13}(1+i+j)=&(1+i+j)(2+i+2j+2k),\\
\beta_{10}(1+i-j)=&(1+i-j)(1+i-2j-2k),\\
\beta_{13}(1+i-j)=&(1+i-j)(2+i-2j-2k),\\
\beta_{10}(1-i+j)=&(1-i+j)(1+i-2j+2k),\\
\beta_{13}(1-i+j)=&(1-i+j)(2+i-2j+2k),\\
\beta_{10}(1-i-j)=&(1-i-j)(1+i+2j-2k),\\
\beta_{13}(1-i-j)=&(1-i-j)(2+i+2j-2k).
\end{align*}
By computation, we have
\begin{scriptsize}
\begin{align*}
(1+i+j)^{-1}\gamma_{r}(1+i+j)=&x_{r}+\frac{-y_{r}-2z_{r}+2t_{r}}{3}i+\frac{-2y_{r}-z_{r}-2t_{r}}{3}j+\frac{-2y_{r}+2z_{r}+t_{r}}{3}k,\\
(1+i-j)^{-1}\gamma_{r}(1+i-j)=&x_{r}+\frac{-y_{r}+2z_{r}-2t_{r}}{3}i+\frac{2y_{r}-z_{r}-2t_{r}}{3}j+\frac{2y_{r}+2z_{r}+t_{r}}{3}k,\\
(1-i+j)^{-1}\gamma_{r}(1-i+j)=&x_{r}+\frac{-y_{r}+2z_{r}+2t_{r}}{3}i+\frac{2y_{r}-z_{r}+2t_{r}}{3}j+\frac{-2y_{r}-2z_{r}+t_{r}}{3}k,\\
(1-i-j)^{-1}\gamma_{r}(1-i-j)=&x_{r}+\frac{-y_{r}-2z_{r}-2t_{r}}{3}i+\frac{-2y_{r}-z_{r}+2t_{r}}{3}j+\frac{2y_{r}-2z_{r}+t_{r}}{3}k,
\end{align*}
\end{scriptsize}
where $r\in\{10,13\}$.

Now, we can naturally get the following lemma.

\begin{lemma}\label{lem3.3}
Suppose $m,n\in\N$.
\begin{enumerate}
\item Given an integer solution $(x_{10}, y_{10},z_{10},t_{10})$ of the system $($$1$-$3$$),$ if either of the following conditions holds:
\begin{itemize}
        \item $y_{10}-z_{10}+t_{10}\equiv0\pmod{3}$,
        \item $y_{10}+z_{10}-t_{10}\equiv0\pmod{3}$,
        \item $y_{10}+z_{10}+t_{10}\equiv0\pmod{3}$,
        \item $y_{10}-z_{10}-t_{10}\equiv0\pmod{3}$,
\end{itemize}
the system $($$1$-$1$-$2$-$2$$)$ has an integer solution.
\item Given an integer solution $(x_{13}, y_{13},z_{13},t_{13})$ of the system $($$2$-$3$$),$   if either of the following conditions holds:
\begin{itemize}
        \item $y_{13}-z_{13}+t_{13}\equiv0\pmod{3}$,
        \item $y_{13}+z_{13}-t_{13}\equiv0\pmod{3}$,
        \item $y_{13}+z_{13}+t_{13}\equiv0\pmod{3}$,
        \item $y_{13}-z_{13}-t_{13}\equiv0\pmod{3}$,
\end{itemize}
the system $($$1$-$2$-$2$-$2$$)$ has an integer solution.
\end{enumerate}
\end{lemma}

In view of Lemma \ref{lem2.1}, we can induce the following lemma.

\begin{lemma}\label{lem3.4}
Suppose $m,n \in \N$.
\begin{enumerate}
\item If $10m-n^{2}\equiv0\ or\ 2 \pmod{3}$, $10m-n^{2}$ is nonnegative and not of the form $4^{r}(8s+7)$ for any $r,s \in \N$,  the system $($$1$-$1$-$2$-$2$$)$ has an integer solution.
\item If $13m-n^{2}\equiv0\ or\ 2 \pmod{3}$, $13m-n^{2}$ is nonnegative and not of the form $4^{r}(8s+7)$ for any $r,s \in \N$,  the system $($$1$-$2$-$2$-$2$$)$ has an integer solution.
\end{enumerate}
\end{lemma}
\noindent{\it Proof}. By Lemma \ref{lem2.1}, either the system ($1$-$1$-$2$-$2$) or the system ($1$-$3$) has an integer solution. If the system ($1$-$1$-$2$-$2$) has an integer solution, we are done. If not, we have an integer solution for the system  ($1$-$3$). Using the notations above, we have
$$\gamma_{10}\beta_{10}=n+A_{10}i+B_{10}j+C_{10}k,$$
where
$$
\begin{cases}
n=x_{10}+3y_{10},\\
A_{10}=3x_{10}-y_{10},\\
B_{10}=-z_{10}-3t_{10},\\
C_{10}=3z_{10}-t_{10}.\\
\end{cases}
$$
We can easily see that
\begin{small}
\begin{align*}
\left(
  \begin{array}{c}  \vspace{1ex}
    x_{10}\\  \vspace{1ex}
    y_{10}\\  \vspace{1ex}
    z_{10}\\  \vspace{1ex}
    t_{10}\\
  \end{array}
\right)
\equiv
  \left(
  \begin{array}{cccc}  \vspace{1ex}
    1 & 0 & 0& 0\\ \vspace{1ex}
    0 & -1 & 0& 0\\ \vspace{1ex}
    0 & 0 & -1&0\\ \vspace{1ex}
    0 &  0  & 0 &-1
  \end{array}
\right)
\left(
  \begin{array}{c}  \vspace{1ex}
   n\\  \vspace{1ex}
   A_{10}\\  \vspace{1ex}
   B_{10}\\  \vspace{1ex}
   C_{10}\\
  \end{array}
\right)  \pmod{3}.
\end{align*}
\end{small}
Hence
\begin{small}
\begin{equation*}
\left(
  \begin{array}{cccc}  \vspace{1ex}
    0 & 1 & -1 & 1\\  \vspace{1ex}
    0 & 1 & 1 & -1\\ \vspace{1ex}
    0 & 1 & 1 & 1\\ \vspace{1ex}
    0 & 1 & -1 & -1\\
  \end{array}
\right)
\left(
  \begin{array}{c}  \vspace{1ex}
    x_{10}\\  \vspace{1ex}
    y_{10}\\  \vspace{1ex}
    z_{10}\\  \vspace{1ex}
    t_{10}\\
  \end{array}
\right)
\equiv
\left(
  \begin{array}{cccc}  \vspace{1ex}
    0 & -1 & 1 & -1\\  \vspace{1ex}
    0 & -1 & -1 & 1\\ \vspace{1ex}
    0 & -1 & -1 & -1\\ \vspace{1ex}
    0 & -1 & 1 & 1\\
  \end{array}
\right)
\left(
  \begin{array}{c}  \vspace{1ex}
   n\\  \vspace{1ex}
   A_{10}\\  \vspace{1ex}
   B_{10}\\  \vspace{1ex}
   C_{10}\\
  \end{array}
\right)  \pmod{3},
\end{equation*}
\end{small}
which means that if
$A_{10}\equiv\pm B_{10}\pm C_{10}\pmod{3}$, the system ($1$-$1$-$2$-$2$) has an integer solution.

On the other hand, $A_{10}^{2}+B_{10}^{2}+C_{10}^{2}=10m-n^{2} \equiv 0\ \t{or}\ 2\pmod{3}$. Note that
$$
\begin{cases}
0\equiv0+0+0\equiv1+1+1\pmod{3},\\
2\equiv0+1+1\pmod{3}.\\
\end{cases}
$$
One can easily check that one of the conditions must hold.
The proof to the second case is similar. \qed
\medskip

{\it Case} 3. ($1$-$4$)$\Rightarrow$($2$-$2$-$3$) and ($2$-$5$)$\Rightarrow$($2$-$3$-$4$)

Let $\beta_{17}=1+4i$, $\beta_{29}=2+5i$, $\gamma_{17}=x_{17}-y_{17}i-z_{17}j-t_{17}k$ and $\gamma_{29}=x_{29}-y_{29}i-z_{29}j-t_{29}k$, where $(x_{17},y_{17},z_{17},t_{17})$ is an integer solution of the system ($1$-$4$) and $(x_{29},y_{29},z_{29},t_{29})$ is an integer solution of the system ($2$-$5$).

One can easily check that
\begin{align*}
\beta_{17}(1+i+j)=&(1-i-j)(-3+2i-2j),\\
\beta_{29}(1+i+j)=&(1-i-j)(-4+3i-2j),\\
\beta_{17}(1+i-j)=&(1-i+j)(-3+2i+2j),\\
\beta_{29}(1+i-j)=&(1-i+j)(-4+3i+2j),\\
\beta_{17}(1-i+j)=&(1+i+j)(3-2i+2k),\\
\beta_{29}(1-i+j)=&(1+i+j)(4-3i+2k),\\
\beta_{17}(1-i-j)=&(1+i-j)(3-2i-2k),\\
\beta_{29}(1-i-j)=&(1+i-j)(4-3i-2k).
\end{align*}
By computation, we have
\begin{scriptsize}
\begin{align*}
(1+i+j)^{-1}\gamma_{r}(1-i-j)=&\frac{-x_{r}-2y_{r}-2z_{r}}{3}+\frac{-2x_{r}-y_{r}+2z_{r}}{3}i+\frac{-2x_{r}+2y_{r}-z_{r}}{3}j+(-t)k,\\
(1+i-j)^{-1}\gamma_{r}(1-i+j)=&\frac{-x_{r}-2y_{r}+2z_{r}}{3}+\frac{-2x_{r}-y_{r}-2z_{r}}{3}i+\frac{2x_{r}-2y_{r}-z_{r}}{3}j+(-t)k,\\
(1-i+j)^{-1}\gamma_{r}(1+i+j)=&\frac{x_{r}+2y_{r}-2t_{r}}{3}+\frac{2x_{r}+y_{r}+2t_{r}}{3}i+(-z)j+\frac{2x_{r}-2y_{r}-t_{r}}{3}k,\\
(1-i-j)^{-1}\gamma_{r}(1+i-j)=&\frac{x_{r}+2y_{r}+2t_{r}}{3}+\frac{2x_{r}+y_{r}-2t_{r}}{3}i+(-z)j+\frac{-2x_{r}+2y_{r}-t_{r}}{3}k,
\end{align*}
\end{scriptsize}
where $r\in\{17,29\}$.

Now, we can naturally get the following lemma.

\begin{lemma}\label{lem3.5}
Suppose $m,n\in \N$.
\begin{enumerate}
\item Given an integer solution $(x_{17}, y_{17},z_{17},t_{17})$ of the system $($$1$-$4$$),$   if either of the following conditions holds:
\begin{itemize}
        \item $x_{17}-y_{17}-z_{17}\equiv0\pmod{3}$,
        \item $x_{17}-y_{17}+z_{17}\equiv0\pmod{3}$,
        \item $x_{17}-y_{17}+t_{17}\equiv0\pmod{3}$,
        \item $x_{17}-y_{17}-t_{17}\equiv0\pmod{3}$,
\end{itemize}
the system $($$2$-$2$-$3$$)$ has an integer solution.
\item Given an integer solution $(x_{29}, y_{29},z_{29},t_{29})$ of the system $($$2$-$5$$),$   if either of the following conditions holds:
\begin{itemize}
        \item $x_{29}-y_{29}-z_{29}\equiv0\pmod{3}$,
        \item $x_{29}-y_{29}+z_{29}\equiv0\pmod{3}$,
        \item $x_{29}-y_{29}+t_{29}\equiv0\pmod{3}$,
        \item $x_{29}-y_{29}-t_{29}\equiv0\pmod{3}$,
\end{itemize}
the system $($$2$-$3$-$4$$)$ has an integer solution.
\end{enumerate}
\end{lemma}
In view of Lemma \ref{lem2.1}, we can induce the following lemma.

\begin{lemma}\label{lem3.6}
Suppose $m,n \in \N$,
\begin{enumerate}
\item If $17m-n^{2}\equiv0\ or\ 2 \pmod{3}$, $17m-n^{2}$ is nonnegative and not of the form $4^{r}(8s+7)$ for any $r,s \in \N$,  the system $($$2$-$2$-$3$$)$ has an integer solution.
\item If $29m-n^{2}\equiv0\ or\ 2 \pmod{3}$, $29m-n^{2}$ is nonnegative and not of the form $4^{r}(8s+7)$ for any $r,s \in \N$,  the system $($$2$-$3$-$4$$)$ has an integer solution.
\end{enumerate}
\end{lemma}

\noindent{\it Proof}. By Lemma \ref{lem2.1}, either the system  $($$2$-$2$-$3$$)$ or the system $($$1$-$4$$)$ has an integer solution. If the system $($$2$-$2$-$3$$)$ has an integer solution,  we are done. If not, we have an integer solution for the system  $($$1$-$4$$)$. Using the notations above, we have
$$\gamma_{17}\beta_{17}=n+A_{17}i+B_{17}j+C_{17}k,$$
where
$$
\begin{cases}
n=x_{17}+4y_{17},\\
A_{17}=4x_{17}-y_{17},\\
B_{17}=-z_{17}-4t_{17},\\
C_{17}=4z_{17}-t_{17}.\\
\end{cases}
$$
We can easily check that
\begin{small}
\begin{equation*}
\left(
  \begin{array}{c}  \vspace{1ex}
    x_{17}\\  \vspace{1ex}
    y_{17}\\  \vspace{1ex}
    z_{17}\\  \vspace{1ex}
    t_{17}\\
  \end{array}
\right)
\equiv
  \left(
  \begin{array}{cccc}  \vspace{1ex}
    -1 & -1 & 0 & 0\\ \vspace{1ex}
    -1 & 1 & 0 & 0\\ \vspace{1ex}
    0 & 0 & 1 & -1\\ \vspace{1ex}
    0 & 0 & 1 &1
  \end{array}
\right)
 \left(
  \begin{array}{c}  \vspace{1ex}
   n\\  \vspace{1ex}
   A_{17}\\  \vspace{1ex}
   B_{17}\\  \vspace{1ex}
   C_{17}\\
  \end{array}
\right)    \pmod{3}.
\end{equation*}
\end{small}
Hence
\begin{small}
\begin{equation*}
\left(
  \begin{array}{cccc}  \vspace{1ex}
    1 & -1 & -1 & 0\\ \vspace{1ex}
    1 & -1 & 1 & 0\\ \vspace{1ex}
    1 & -1 & 0 & 1\\ \vspace{1ex}
    1 & -1 & 0 & -1\\
  \end{array}
\right)
 \left(
  \begin{array}{c}  \vspace{1ex}
    x_{17}\\  \vspace{1ex}
    y_{17}\\  \vspace{1ex}
    z_{17}\\  \vspace{1ex}
    t_{17}\\
  \end{array}
\right)
\equiv
\left(
  \begin{array}{cccc} \vspace{1ex}
    0 & 1 & -1 & 1\\ \vspace{1ex}
    0 & 1 & 1 & -1\\ \vspace{1ex}
    0 & 1 & 1 & 1\\ \vspace{1ex}
    0 & 1 & -1 & -1\\
  \end{array}
\right)
\left(
  \begin{array}{c}  \vspace{1ex}
   n\\  \vspace{1ex}
   A_{17}\\  \vspace{1ex}
   B_{17}\\  \vspace{1ex}
   C_{17}\\
  \end{array}
\right)    \pmod{3},
\end{equation*}
\end{small}
which means that if $A_{17}\equiv\pm B_{17}\pm C_{17}\pmod{3}$, the system $($$2$-$2$-$3$$)$ has an integer solution.

On the other hand, $A_{17}^{2}+B_{17}^{2}+C_{17}^{2}=17m-n^{2} \equiv 0\ \t{or}\ 2\pmod{3}$ . Note that
$$
\begin{cases}
0\equiv0+0+0\equiv1+1+1\pmod{3},\\
2\equiv0+1+1\pmod{3}.\\
\end{cases}
$$
One can easily check that one of the conditions must hold.
The proof to the second case is similar. \qed
\medskip

{\it Case} 4. $($$a$-$a$-$a$-$b$$)$ $\Rightarrow$ $($$\frac{a+4b}{5}$-$\frac{7a-2b}{5}$-$\frac{3a+2b}{5}$-$\frac{4a+b}{5}$$)$, where $a,b$ are distinct integers satisfying that  $a\equiv b\pmod{5}$ and $ab$ is not divided by $5$.

Let $\eta=a+ai+aj+bk$ and $\gamma=x_{0}-y_{0}i-z_{0}j-t_{0}k$, where $(x_{0},y_{0},z_{0},t_{0})$ is an integer solution of the system $($$a$-$a$-$a$-$b$$)$.

One can easily check that
\begin{align*}
(1+2j)^{-1}\eta(1+2i)=\frac{a+4b}{5}+\frac{7a-2b}{5}i+\frac{3a+2b}{5}j+\frac{4a+b}{5}k,\\
(1-2j)^{-1}\eta(1-2i)=\frac{a+4b}{5}+\frac{3a+2b}{5}i+\frac{7a-2b}{5}j+\frac{4a+b}{5}k,\\
(1+2k)^{-1}\eta(1+2j)=\frac{3a+2b}{5}+\frac{7a-2b}{5}i+\frac{a+4b}{5}j+\frac{4a+b}{5}k,\\
(1-2k)^{-1}\eta(1-2j)=\frac{7a-2b}{5}+\frac{3a+2b}{5}i+\frac{a+4b}{5}j+\frac{4a+b}{5}k,\\
(1+2i)^{-1}\eta(1+2k)=\frac{7a-2b}{5}+\frac{a+4b}{5}i+\frac{3a+2b}{5}j+\frac{4a+b}{5}k,\\
(1-2i)^{-1}\eta(1-2k)=\frac{3a+2b}{5}+\frac{a+4b}{5}i+\frac{7a-2b}{5}j+\frac{4a+b}{5}k.
\end{align*}
By computation, we have
\begin{align*}
(1+2i)^{-1}\gamma(1+2j)=&\frac{x_{0}-2y_{0}+2z_{0}+4t_{0}}{5}+\frac{-2x_{0}-y_{0}-4z_{0}+2t_{0}}{5}i\\
&+\frac{2x_{0}-4y_{0}-z_{0}-2t_{0}}{5}j+\frac{-4x_{0}-2y_{0}+2z_{0}-t_{0}}{5}k,\\
(1-2i)^{-1}\gamma(1-2j)=&\frac{x_{0}+2y_{0}-2z_{0}+4t_{0}}{5}+\frac{2x_{0}-y_{0}-4z_{0}-2t_{0}}{5}i\\
&+\frac{-2x_{0}-4y_{0}-z_{0}+2t_{0}}{5}j+\frac{-4x_{0}+2y_{0}-2z_{0}-t_{0}}{5}k,\\
(1+2j)^{-1}\gamma(1+2k)=&\frac{x_{0}+4y_{0}-2z_{0}+2t_{0}}{5}+\frac{-4x_{0}-y_{0}-2z_{0}+2t_{0}}{5}i\\
&+\frac{-2x_{0}+2y_{0}-z_{0}-4t_{0}}{5}j+\frac{2x_{0}-2y_{0}-4z_{0}-t_{0}}{5}k,\\
(1-2j)^{-1}\gamma(1-2k)=&\frac{x_{0}+4y_{0}+2z_{0}-2t_{0}}{5}+\frac{-4x_{0}-y_{0}+2z_{0}-2t_{0}}{5}i\\
&+\frac{2x_{0}-2y_{0}-z_{0}-4t_{0}}{5}j+\frac{-2x_{0}+2y_{0}-4z_{0}-t_{0}}{5}k,\\
(1+2k)^{-1}\gamma(1+2i)=&\frac{x_{0}+2y_{0}+4z_{0}-2t_{0}}{5}+\frac{2x_{0}-y_{0}-2z_{0}-4t_{0}}{5}i\\
&+\frac{-4x_{0}+2y_{0}-z_{0}-2t_{0}}{5}j+\frac{-2x_{0}-4y_{0}+2z_{0}-t_{0}}{5}k,\\
(1-2k)^{-1}\gamma(1-2i)=&\frac{x_{0}-2y_{0}+4z_{0}+2t_{0}}{5}+\frac{-2x_{0}-y_{0}+2z_{0}-4t_{0}}{5}i\\
&+\frac{-4x_{0}-2y_{0}-z_{0}+2t_{0}}{5}j+\frac{2x_{0}-4y_{0}-2z_{0}-t_{0}}{5}k.
\end{align*}

Now, we can naturally get the following lemma.

\begin{lemma}\label{lem3.7}
Given an integer solution $(x_{0}, y_{0},z_{0},t_{0})$ of the system $($$a$-$a$-$a$-$b$$)$, if either of the following conditions holds:
\begin{itemize}
        \item $x_{0}-2y_{0}+2z_{0}-t_{0}\equiv0\pmod{5}$,
        \item $x_{0}+2y_{0}-2z_{0}-t_{0}\equiv0\pmod{5}$,
        \item $x_{0}-y_{0}-2z_{0}+2t_{0}\equiv0\pmod{5}$,
        \item $x_{0}-y_{0}+2z_{0}-2t_{0}\equiv0\pmod{5}$,
        \item $x_{0}+2y_{0}-z_{0}-2t_{0}\equiv0\pmod{5}$,
        \item $x_{0}-2y_{0}-z_{0}+2t_{0}\equiv0\pmod{5}$,
\end{itemize}
the system $($$\frac{a+4b}{5}$-$\frac{7a-2b}{5}$-$\frac{3a+2b}{5}$-$\frac{4a+b}{5}$$)$ has an integer solution.
\end{lemma}

In view of Lemma \ref{lem2.1}, we can induce the following lemma.

\begin{lemma}\label{lem3.8}
Let $a$,$b$ be two distinct integers satisfying $a\equiv b\pmod{5}$ and $ab$ is not divided by $5$. Suppose that $l=3a^{2}+b^{2}$ has only two ways to be represented as sum of four squares. Assume that  $m,n\in \N$ satisfy that $lm-n^{2}$ is nonnegative and not of the form of $4^{r}(8s+7)$ for any $r,s \in \N$. When $lm-n^{2}\equiv 0,1\ or\ 4\pmod{5}$, the system $($$\frac{a+4b}{5}$-$\frac{7a-2b}{5}$-$\frac{3a+2b}{5}$-$\frac{4a+b}{5}$$)$  has an integer solution.
\end{lemma}

\noindent{\it Proof}. By Lemma \ref{lem2.1}, either the system $($$\frac{a+4b}{5}$-$\frac{7a-2b}{5}$-$\frac{3a+2b}{5}$-$\frac{4a+b}{5}$$)$  or the system $($$a$-$a$-$a$-$b$$)$ has an integer solution. If the system $($$\frac{a+4b}{5}$-$\frac{7a-2b}{5}$-$\frac{3a+2b}{5}$-$\frac{4a+b}{5}$$)$  has an integer solution, we are done. If not, we have an integer solution of the system  $($$a$-$a$-$a$-$b$$)$. Using the notations above, we have
$$\gamma(a+ai+aj+bk)=n+Ai+Bj+Ck,$$
where
$$
\begin{cases}
n=ax_{0}+ay_{0}+az_{0}+bt_{0},\\
A=ax_{0}-ay_{0}-bz_{0}+at_{0},\\
B=ax_{0}+by_{0}-az_{0}-at_{0},\\
C=bx_{0}-ay_{0}+az_{0}-at_{0}.\\
\end{cases}
$$
By solving the equations and taking the coefficients module $5$, we have
\begin{small}
\begin{equation*}
\left(
  \begin{array}{c}  \vspace{1ex}
    x_{0}\\  \vspace{1ex}
    y_{0}\\  \vspace{1ex}
    z_{0}\\  \vspace{1ex}
    t_{0}\\
  \end{array}
\right)
\equiv
a^{-1}\left(
  \begin{array}{cccc} \vspace{1ex}
    -1 & -1 & -1 & -1\\  \vspace{1ex}
    -1 & 1 & -1 & 1\\  \vspace{1ex}
    -1 & 1 & 1 & -1\\  \vspace{1ex}
    -1 & -1 & 1 & 1
\end{array}
\right)
  \left(
  \begin{array}{c}  \vspace{1ex}
   n\\  \vspace{1ex}
   A\\  \vspace{1ex}
   B\\  \vspace{1ex}
   C\\
  \end{array}
\right)
\pmod{5}.
\end{equation*}
\end{small}
Hence
\begin{small}
\begin{equation*}
\left(
  \begin{array}{cccc}   \vspace{1ex}
    1 & -2 & 2 & -1\\  \vspace{1ex}
    1 & 2 & -2 & -1\\ \vspace{1ex}
    1 & -1 & -2 & 2\\ \vspace{1ex}
    1 & -1 & 2 & -2\\ \vspace{1ex}
    1 & 2 & -1 & -2\\ \vspace{1ex}
    1 & -2 & -1 & 2\\
  \end{array}
\right)
\left(
  \begin{array}{c}  \vspace{1ex}
    x_{0}\\  \vspace{1ex}
    y_{0}\\  \vspace{1ex}
    z_{0}\\  \vspace{1ex}
    t_{0}\\
  \end{array}
\right)
\equiv
a^{-1}
\left(
  \begin{array}{cccc} \vspace{1ex}
    0 & 0 & 2 & -1\\  \vspace{1ex}
    0 & 0 & -1 & 2\\ \vspace{1ex}
    0 & -1 & 0 & 2\\ \vspace{1ex}
    0 & 2 & 0 & -1\\ \vspace{1ex}
    0 & 2 & -1 & 0\\ \vspace{1ex}
    0 & -1 & 2 & 0\\
  \end{array}
\right)
\left(
  \begin{array}{c}  \vspace{1ex}
   n\\  \vspace{1ex}
   A\\  \vspace{1ex}
   B\\  \vspace{1ex}
   C\\
  \end{array}
\right)
\pmod{5},
\end{equation*}
\end{small}
which means that if either of the following conditions holds:
\begin{itemize}
  \item $A\equiv \pm 2B \pmod{5}$,
  \item $A\equiv \pm 2C \pmod{5}$,
  \item $B\equiv \pm 2C \pmod{5}$,
\end{itemize}
the system $($$\frac{a+4b}{5}$-$\frac{7a-2b}{5}$-$\frac{3a+2b}{5}$-$\frac{4a+b}{5}$$)$ has an integer solution.

On the other hand, $A^{2}+B^{2}+C^{2}=lm-n^{2} \equiv 0,1\ \t{or} \ 4\pmod{5}$. Note that

$$
\begin{cases}
0\equiv0+0+0\equiv0+1+4\pmod{5},\\
1\equiv0+0+1\equiv1+1+4\pmod{5},\\
4\equiv0+0+4\equiv1+4+4\pmod{5}.
\end{cases}
$$
One can easily check that one of $A^{2}+B^{2}$, $A^{2}+C^{2}$ and $B^{2}+C^{2}$ must be divided by $5$, hence one of the conditions must hold. This completes the proof. \qed
\medskip

\begin{remark}
	In fact, we can replace $(a,b)$ with $(1,-4)$, $(1,6)$, $(2,-3)$ and $(3,-2)$, which gives the existance of integer solutions of the system $(1$-$3$-$3)$, $(1$-$2$-$3$-$5)$, $(1$-$2$-$4)$ and $(1$-$1$-$2$-$5)$ under certain conditions. However, we need to do some more work with the case $($$1$-$1$-$1$-$6$$)$$\Rightarrow$($1$-$2$-$3$-$5$) to prove Theorem\ref{Th1.4}.
\end{remark}

{\it Case} 5. $($$1$-$1$-$1$-$6$$)$$\Rightarrow$($1$-$2$-$3$-$5$)

Let $\beta_{39}=1+i+j+6k$ and $\gamma_{39}=x_{39}-y_{39}i-z_{39}j-t_{39}k$, where $(x_{39},y_{39},z_{39},t_{39})$ is an integer solution of the system ($1$-$1$-$1$-$6$).

One can easily check that
\begin{align*}
\beta_{39}(1+i+j)=&(1+i+j)(1-3i+5j-2k),\\
\end{align*}
By computation, we have
\begin{small}
\begin{align*}
(1+i+j)^{-1}\gamma_{39}(1+i+j)=&x_{39}+\frac{-y_{39}-2z_{39}+2t_{39}}{3}i\\
&+\frac{-2y_{39}-z_{39}-2t_{39}}{3}j+\frac{-2y_{39}+2z_{39}+t_{39}}{3}k,\\
\end{align*}
\end{small}
Note that if $(x_{39}, y_{39},z_{39},t_{39})$ is a solution of the system ($1$-$1$-$1$-$6$), $(x_{39}, z_{39},y_{39},t_{39})$, $(y_{39}, x_{39},z_{39},t_{39})$, $(y_{39}, z_{39},x_{39},t_{39})$, $(z_{39}, x_{39},y_{39},t_{39})$ and $(z_{39}, y_{39},x_{39},t_{39})$ are also solutions of the system  ($1$-$1$-$1$-$6$). Hence we can naturally get the following lemma.
\begin{lemma}\label{lem3.9}
 Given an integer solution $(x_{39}, y_{39},z_{39},t_{39})$ for the system $($$1$-$1$-$1$-$6$$)$,  if either of the following conditions holds:
\begin{itemize}
        \item $t_{39}\equiv z_{39}-y_{39}\pmod{3}$,
        \item $t_{39}\equiv x_{39}-z_{39}\pmod{3}$,
        \item $t_{39}\equiv y_{39}-x_{39}\pmod{3}$,
        \item $t_{39}\equiv y_{39}-z_{39}\pmod{3}$,
        \item $t_{39}\equiv x_{39}-y_{39}\pmod{3}$,
        \item $t_{39}\equiv z_{39}-x_{39}\pmod{3}$,
\end{itemize}
the system $($$1$-$2$-$3$-$5$$)$ has an integer solution.
\end{lemma}

In view of Lemma \ref{lem2.1}, we can induce the following lemma.

\begin{lemma}\label{lem3.10} Suppose $m,n \in \N$ and $n\not\equiv 0\pmod{3}$. If $39m-n^{2}$ is nonnegative and not of the form $4^{r}(8s+7)$ for any $r,s \in \N$, the system $($$1$-$2$-$3$-$5$$)$ has an integer solution.
\end{lemma}

\noindent{\it Proof}. By Lemma \ref{lem2.1}, either the system ($1$-$2$-$3$-$5$) or the system ($1$-$1$-$1$-$6$) has an integer solution. If the system ($1$-$2$-$3$-$5$) has an integer solution, we are done. If not, we have an integer solution for the system  ($1$-$1$-$1$-$6$), say $(x_{39}, y_{39},z_{39},t_{39})$. Suppose that none of the six conditions in Lemma \ref{lem3.9} holds, let us see what happens in the following cases:
\begin{itemize}
\item $t_{39}\equiv0\pmod{3}$. It is easy to see that $x_{39}, y_{39}$ and $z_{39}$ are pairwise distinct module $3$. Thus $x_{39}+y_{39}+z_{39}\equiv0\pmod{3}$. However, $x_{39}+y_{39}+z_{39}\equiv x_{39}+y_{39}+z_{39}+6t_{39} \equiv n \not\equiv 0\pmod{3}$, which is a contradiction.
\item $t_{39}\not\equiv 0\pmod{3}$. Since $x_{39}-y_{39}\not\equiv \pm t_{39}\pmod{3}$, we have $x_{39}\equiv y_{39}\pmod{3}$. Similarly, we have $x_{39}\equiv y_{39}\equiv z_{39}\pmod{3}$. Hence, $x_{39}+y_{39}+z_{39}\equiv0\pmod{3}$,  which is again a contradiction.
\end{itemize}
\qed
\medskip

\section{Proof of the theroems}
In this section we set $l=a^{2}+b^{2}+c^{2}+d^{2}$, where $(a,b,c,d)$ is among the quadruples mentioned in the theorems. For $n\in\N$, $\ord_{2}(n)$ denotes the integer satisfying that $2^{\ord_{2}(n)}\mid n$ and $2^{\ord_{2}(n)+1}\nmid n$.

\medskip
\noindent{\it Proof of Theorem} 1.1.
For $m\in\N$, define $C_{m}$ as $\{n \in \N:\ m-n^{6}\geq 0 \ \t{and} \ m-n^{6}\ \t{is not of the}$ $ \t{form}\ 4^{r}(8s+7)\ \t{for any}\ r,s \in \N \}$. We will discuss the elements of $C_{m}$.
\begin{itemize}
\item If $m\equiv 1,2,3,5,6\pmod{8}$, all nonnegative even integers less than $\sqrt[6]{m}$ are in $C_{m}$ by simple congruence argument.
\item If $m\equiv 2,3,4,6,7\pmod{8}$, all nonnegative odd integers less than $\sqrt[6]{m}$ are in $C_{m}$ by simple congruence argument.
\item If $\ord_{2}(m)=3\ or\ 5$, $\ord_{2}(m-n^{6})$ must be odd for any even integer $n$, hence all even nonnegative integers less than $\sqrt[6]{m}$ are in $C_{m}$.
\item If $\ord_{2}(m)=4$, for nonnegative even integer $n$, $m-n^{6}=2^{4}(m/16-2^{2}(n/2)^{6})$. After simple congruence argument we can get that if $m/16\equiv 1,3,5\pmod{8}$, $\{ n\in\N:m-n^{6}\geq 0 \ \t{and} \ n\equiv0\pmod{4}\}\subseteq C_{m}$. If $m/16\equiv 1,5,7\pmod{8}$, $\{n\in\N: m-n^{6}\geq 0 \ \t{and} \ n\equiv2\pmod{4}\}\subseteq C_{m}$.
\item If $\ord_{2}(m)=6$, for nonnegative even integer $n$, $m-n^{6}=2^{6}(m/64-(n/2)^{6})$. After simple congruence argument we can get that if $m/64\equiv 1,3,5\pmod{8}$, $\{ n\in\N:m-n^{6}\geq 0 \ \t{and} \ n\equiv0\pmod{4}\}\subseteq C_{m}$. If $m/64\equiv 3,5,7\pmod{8}$, $\{n\in\N: m-n^{6}\geq 0 \ \t{and} \ n\equiv2\pmod{4}\}\subseteq C_{m}$.

\end{itemize}
Therefore, when $m$ is not divided by 64 and $\sqrt[6]{lm}\geq 10$, there is $n\in C_{lm}$ such that $n^{3}$ satisfies the conditions in Lemma \ref{lem3.2}, Lemma \ref{lem3.4}, Lemma \ref{lem3.6} and Lemma \ref{lem3.8}. For $m$ less than $10^{6}/l$ we can easily check by computer. When $m$ is divided by $64$, it suffices to consider $m/64$.
\qed
\medskip

\medskip
\noindent{\it Proof of Theorem} 1.2.
For $m\in\N$, define $P_{m}$ as $\{k \in \N:\ m-4^{k}\geq 0 \ \t{and} \ m-4^{k}\ \t{is}\ \t{not}\ \t{of}\ \t{the}\ \t{form}\ 4^{r}(8s+7)\ \t{for any}\ r,s \in \N \}$. We will discuss the elements of $P_{m}$.

\begin{itemize}
\item If $m\equiv 1\pmod{4}$, any integer $k$ satisfying $1\leq k \leq \log_4m$ is in $P_{m}$ by simple congruence argument.
\item If $m\equiv 2\pmod{4}$, any integer $k$ satisfying $0\leq k \leq \log_4m$ is in $P_{m}$ by simple congruence argument.
\item If $m\equiv 3\pmod{8}$, any integer $k$ satisfying $0\leq k \leq \log_4m$ and $k\neq1$ is in $P_{m}$ by simple congruence argument.
\item If $m\equiv 7\pmod{8}$, $\{0,1\}\subseteq P_{m}$ by simple congruence argument.
\item If $\ord_{2}(m)=2$, it is easy to see that $0\in P_{m}$. Since $m-4^{k}=4(m/4-4^{k-1})$, one can easily check that $1\in P_{m}$ when $m/4 \not\equiv 5\pmod{8}$ and any integer $k$ satisfying $3\leq k\leq \log_4m$ is in $P_{m}$ when $m/4 \not\equiv 7\pmod{8}$.
\end{itemize}
Therefore, when $m$ is not divided by $4$ and $m\geq 64/l$, there is $k \in P_{lm}$ such that $2^{k}$ satisfies the conditions in Lemma \ref{lem3.2} and Lemma \ref{lem3.8}. For $m$ less than $64/l$, we can easily check by computer. For $m$ divided
 by $4$, it suffices to consider $m/4$.  \qed
\medskip

\medskip
\noindent{\it Proof of Theorem} 1.3.
For $m\in\N$, define $S_{m}$ as $\{n \in \N:\ m-n^{4}\geq 0 \ \t{and} \ m-n^{4}\ \t{is not of the}$ $ \t{form}\ 4^{r}(8s+7)\ \t{for any}\ r,s \in \N \}$. We will discuss the elements of $S_{m}$.
\begin{itemize}
\item If $m\equiv 1,2,3,5,6\pmod{8}$, all nonnegative even integers less than $\sqrt[4]{m}$ are in $S_{m}$ by simple congruence argument.
\item If $m\equiv 2,3,4,6,7\pmod{8}$, all nonnegative odd integers less than $\sqrt[4]{m}$ are in $S_{m}$ by simple congruence argument.
\item If $\ord_{2}(m)=3$, then $\ord_{2}(m-n^{4})$ must be odd for any even integer $n$, hence all even nonnegative integers less than $\sqrt[4]{m}$ are in $S_{m}$.
\item If $\ord_{2}(m)=4$, for nonnegative even integer $n$, $m-n^{4}=2^{4}(m/16-(n/2)^{4})$. After simple congruence argument we can get that if $m/16\equiv 1,3,5\pmod{8}$, $\{ n\in\N:m-n^{4}\geq 0 \ \t{and} \ n\equiv0\pmod{4}\}\subseteq S_{m}$. If $m/16\equiv 1,5,7\pmod{8}$, $\{n\in\N: m-n^{4}\geq 0 \ \t{and} \ n\equiv2\pmod{4}\}\subseteq S_{m}$.
\end{itemize}
Therefore, if $\sqrt[4]{lm}\geq 6$ and $m\not\equiv0\pmod{16}$, there is $n\in S_{lm}$ such that $n^{2}$ satisfies the conditions in Lemma \ref{lem3.4} and Lemma \ref{lem3.6}. For $m$ less than $6^{4}/l$ we can easily check by computer. For $m$ divided by 16, it suffice to consider $m/16$.  \qed
\medskip

\noindent{\it Proof to Theorem} 1.4. When $m>3.74\times10^{9}>(4/(\sqrt[4]{39}-\sqrt[4]{38}))^{4}$ and $m\not\equiv0\pmod{16}$, we can find an integer $n \in S_{39m}\bigcap [\sqrt[4]{38m},\sqrt[4]{39m}]$ satisfying that $n\not\equiv0\pmod{3}$, hence the system ($1$-$2$-$3$-$5$) has an integer solution, say ($x_{39},y_{39},z_{39},t_{39}$). In view of Lemma 2.1 in \cite{S4}, ($x_{39},y_{39},z_{39},t_{39}$) is in fact a natural solution. To prove the second claim, it suffice to note that $7.67\times10^{9}>(6/(\sqrt[4]{29}-\sqrt[4]{28}))^{4}$.
\qed

\begin{acknowledgment}
The authors would like to thank Prof. Zhi-Wei Sun and Dr. Chen Wang for their helpful comments. This work is supported by the National Natural Science Foundation of China (Grant No. 11971222).
\end{acknowledgment}


\begin{thebibliography}{KMP}
\bibitem{D1} L. E. Dickson, {\it Quaternary quadratic forms representing all integers}, Amer. J. Math. {\bf 49} (1927), 39--56.
\bibitem{D} L. E. Dickson, {\it Modern Elementary Theory of Numbers}, University of Chicago Press, Chicago, 1939.
\bibitem{G} E. Grosswald, {\it Representation of Integers as Sums of Squares}, Springer, New York, 1985.
\bibitem{MT} A. Machiavelo and N. Tsopanidis, {\it Zhi-Wei Sun's 1-3-5 Conjecture and Variations}, preprint, arXiv:2003.02592, 2020.
\bibitem{P} G. Pall, {\it On the Arithmetic of Quaternions}, Trans. Amer. Math. Soc. {\bf 47} (1940), 487--500.
\bibitem{R} S. Ramanujan, {\it On the expression of a number in the from $ax^{2} + by^{2} + cz^{2} + dv^{2}$}, Proc. Cambridge Philos. Soc. {\bf 19} (1917), 11--21.
\bibitem{SS}Y.-C. Sun and Z.-W. Sun, {\it Some variants of Lagrange's four squares theorem}, Acta Arith. {\bf 183} (2018), 339--356.
\bibitem{S1}Z.-W. Sun, {\it On universal sums of polygonal numbers}, Sci. China Math. {\bf 58} (2015), 1367--1396.
\bibitem{S2}Z.-W. Sun, {\it Refining Lagrange's four-square theorem}, J. Number Theory {\bf 175} (2017), 169--190.
\bibitem{S3}Z.-W. Sun, {\it Restricted sums of four squares}, Int. J. Number Theory {\bf 15} (2019), 1863--1893.
\bibitem{S4}Z.-W. Sun, {\it Sums of four squares with certain restrictions}, preprint, arXiv:2010.05775, 2020.
\bibitem{W} H.-L. Wu and Z.-W. Sun, {\it On the 1-3-5 conjecture and related topics}, Acta Arith. {\bf 193} (2020), 253--268.
\end{thebibliography}
\end{document}